\documentclass[12pt]{article}
\usepackage{amsfonts,amsmath,epsfig}
\usepackage{mitpress}

\newdimen\dummy
\dummy=\oddsidemargin
\begin{document}

\begin{center}
\textbf{A twisted group algebra structure for an algebra obtained by the
Cayley-Dickson process}%
\begin{equation*}
\end{equation*}

\bigskip Cristina FLAUT and Remus BOBOESCU%
\begin{equation*}
\end{equation*}
\end{center}

\textbf{Abstract.} {\small Starting from some ideas given by Bales in [Ba;
09], in this paper we present an algorithm for computing the elements of the
basis in an algebra obtained by the Cayley-Dickson process. As a consequence
of this result, we prove that an algebra obtained by the Cayley-Dickson
process is a twisted group algebra for the group }$G=\mathbb{Z}%
_{2}^{n},n=2^{t}${\small , }$t\in \mathbb{N}${\small , over a field }$K$%
{\small , with char}$K=0${\small . In the last section, we give some
properties and applications of the quaternion nonassociative
algebras.\medskip }%
\begin{equation*}
\end{equation*}

\textbf{2010} \textbf{AMS Classification}: 17A35


\begin{equation*}
\end{equation*}

\textbf{1. Introduction}%
\begin{equation*}
\end{equation*}

In the following, we consider $K$ a commutative field with $charK\neq 2$ and 
$\mathcal{E}$ an algebra over the field \ $K$. An algebra $\mathcal{E}$ is
called \textit{unitary} if this algebra contains an identity element with
respect to the algebra's multiplication.

The set 
\begin{equation*}
\mathcal{N}\left( \mathcal{E}\right) =\{x\in \mathcal{E}~/~\left(
x,a,b\right) =\left( a,x,b\right) =\left( a,b,x\right) =0\text{, for all }%
a,b\in \mathcal{E}\},
\end{equation*}%
is called \textit{the nucleus} of the algebra $\mathcal{E}$.

An algebra $\mathcal{E}$ is called \textit{alternative} if $x^{2}y=x\left(
xy\right) $ and $xy^{2}=\left( xy\right) y,$ for all $x,y\in \mathcal{E},$ 
\textit{\ flexible} if $x\left( yx\right) =\left( xy\right) x=xyx,$ for all $%
x,y\in \mathcal{E}$ and \textit{power associative} if the subalgebra $<x>$
of $\mathcal{E}$, generated by any element $x\in \mathcal{E}$, is
associative. $\ $Each alternative algebra is$\ $a\ flexible algebra and a
power associative algebra. A unitary algebra $\mathcal{E}\neq K$ such that
the following relation $x^{2}+\alpha _{x}x+\beta _{x}=0$ is true for each $%
x\in \mathcal{E},$ with $\alpha _{x},\beta _{x}\in K,$ is called a \textit{%
quadratic algebra}. A finite-dimensional algebra $\mathcal{E}$ is \textit{a
division} algebra if and only if $\mathcal{E}$ does not contain zero
divisors. (See [Sc; 66])\medskip

In the following, we briefly present the \textit{Cayley-Dickson process} and
the properties of the obtained algebras. (see [Sc; 66] and [Sc; 54]).

We consider $\mathcal{E},$ a finite dimensional unitary algebra over a field 
$\ K,$ with a \textit{scalar} \textit{involution} $\,$%
\begin{equation*}
\,\,\,\overline{\phantom{x}}:\mathcal{E}\rightarrow \mathcal{E},a\rightarrow 
\overline{a},
\end{equation*}%
$\,\,$ which it is a linear map with the following properties$\,\,\,$%
\begin{equation*}
\overline{ab}=\overline{b}\overline{a},\,\overline{\overline{a}}=a,
\end{equation*}%
$\,\,$and 
\begin{equation*}
a+\overline{a},a\overline{a}\in K\cdot 1\ \text{for all }a,b\in \mathcal{E}.%
\text{ }
\end{equation*}%
An element $\,\overline{a}$ is called the \textit{conjugate} of the element $%
a$. The linear form$\,\,$%
\begin{equation*}
\,\,\mathbf{t}:A\rightarrow K\,,\,\,\mathbf{t}\left( a\right) =a+\overline{a}
\end{equation*}%
and the quadratic form 
\begin{equation*}
\mathbf{n}:A\rightarrow K,\,\,\mathbf{n}\left( a\right) =a\overline{a}\ 
\end{equation*}%
are called the \textit{trace} and the \textit{norm \ }of \ the element $a$,
respectively. From here, it results that an algebra $\mathcal{E}$ with a
scalar involution is quadratic. $\,$

We consider$\,\,\,\gamma \in K$ \thinspace\ a fixed non-zero element. We
define the following algebra multiplication on the vector space 
\begin{equation}
\mathcal{E}\oplus \mathcal{E}:\left( a_{1},a_{2}\right) \left(
b_{1},b_{2}\right) =\left( a_{1}b_{1}+\gamma \overline{b}_{2}a_{2},a_{2}%
\overline{b_{1}}+b_{2}a_{1}\right) .  \tag{1.1}
\end{equation}%
\newline
The obtained algebra structure over $\mathcal{E}\oplus \mathcal{E},$ denoted
by $\left( \mathcal{E},\gamma \right) $ is called the \textit{algebra
obtained from }$\mathcal{A}$\textit{\ by the Cayley-Dickson process.} $\,$We
have $\dim \left( \mathcal{E},\gamma \right) =2\dim \mathcal{E}$.

Let $x\in \left( \mathcal{E},\gamma \right) $, $x=\left( a_{1},a_{2}\right) $%
. The map 
\begin{equation*}
\,\,\,\overline{\phantom{x}}:\left( \mathcal{E},\gamma \right) \rightarrow
\left( \mathcal{E},\gamma \right) \,,\,\,x\rightarrow \bar{x}\,=\left( 
\overline{a}_{1},-a_{2}\right) ,
\end{equation*}%
\newline
is a scalar involution of the algebra $\left( \mathcal{E},\gamma \right) $,
extending the involution $\overline{\phantom{x}}\,\,\,$of the algebra $%
\mathcal{E}$. We have that 
\begin{equation*}
\,\mathbf{t}\left( x\right) =\mathbf{t}(a_{1})
\end{equation*}%
and$\,\,\,$ 
\begin{equation*}
\mathbf{n}\left( x\right) =\mathbf{n}\left( a_{1}\right) -\gamma \mathbf{n}%
(a_{2})
\end{equation*}%
are the \textit{trace} and the \textit{norm} of the element $x\in $ $\left( 
\mathcal{E},\gamma \right) $, respectively.\thinspace $\,$

\thinspace If we consider $\mathcal{E}=K$ \thinspace and we apply this
process $t$ times, $t\geq 1$,$\,\,$we obtain an algebra over $K$,$\,\,$%
\begin{equation}
\mathcal{E}_{t}=\left( \frac{\gamma _{1},...,\gamma _{t}}{K}\right) \text{.}
\tag{1.2. }
\end{equation}

Using induction in this algebra, the set $\{1,f_{1},...,f_{n-1}\},n=2^{t}$,
generates a basis with the properties:%
\begin{equation}
f_{i}^{2}=\gamma _{i}1,\,\,_{i}\in K,\gamma _{i}\neq 0,\,\,i\in \{1,2..,n-1\}
\tag{1.3.}
\end{equation}%
and \ 
\begin{equation}
f_{i}f_{j}=-f_{j}f_{i}=\beta _{ij}f_{k},\,\,\beta _{ij}\in K,\,\,\beta
_{ij}\neq 0,i\neq j,i,j\in \{1,2,...,n-1\},  \tag{1.4.}
\end{equation}%
$\ \ \beta _{ij}$ and $f_{k}$ being uniquely determined by $f_{i}$ and $f_{j}
$.(See [Sc; 54]).

For $t=2,$ we obtain the generalized quaternion algebras and for $t=3$, we
obtain the generalized octonion algebras.%
\begin{equation*}
\end{equation*}

\textbf{2. Twisted group algebra structure for the algebra} $\mathcal{E}_{t}$
\begin{equation*}
\end{equation*}

In [Ba; 09], was presented the twist tree for an algebra obtained by the
Cayley-Dickson process in the real case and in the situation when$~\gamma
_{i}=-1$, for all $i\in \{2,...,n\}$. As a consequence, an algorithm for
computing two elements of the basis in this particular case was provided. In
the following, we will prove that the algebra $\mathcal{E}_{t}$ has a
twisted group algebra structure. Moreover, we give an algorithm which allow
us to compute two elements of the basis, in the general case of the algebra $%
\mathcal{E}_{t}$. In this way the calculations become more easier in higher
dimension of the algebra $\mathcal{E}_{t}$.\medskip 

\textbf{Definition} \textbf{2.1.} ( [Re; 71]) Let \thinspace $\left( G,\cdot
\right) $ be a finite group and $K$ be a field. A \textit{twisted group
algebra} for the group $G$ over the field $K$ is an algebra over the field $%
K $ with a basis $\{a_{g},g\in G\}$ such that 
\begin{equation*}
a_{g}a_{h}=f\left( g,h\right) a_{g\cdot h},\text{ where }g,h\in G,f\left(
g,h\right) \in K,f\left( g,h\right) \neq 0.
\end{equation*}

\textbf{Remark 2.2.} If we consider $\gamma _{1}=...=\gamma _{t}=-1$ and $K=%
\mathbb{R}$, in [Ba; 09], was described \ how the basis vectors can be
multiplied in the algebra $\mathcal{E}_{t}$, with $\dim \mathcal{E}%
_{t}=2^{t}=n$. He used the binary decomposition for the subscript indices.

If $\ f_{p},f_{q}$ are two vectors of the basis $B$ with $p,q$ representing
the binary decomposition for the indices of the vectors, that means $p,q$
are in $\mathbb{Z}_{2}^{n}$, we have that $f_{p}f_{q}=\alpha _{t}\left(
p,q\right) f_{p\otimes q}=f_{pq},$ where:

i) $p\otimes q$ are the sum of \ $p$ and $q$ in the group $\mathbb{Z}
_{2}^{n} $ or, more precisely, the "\textit{exclusive or}" for the binary
numbers $p$ and $q;$

ii) $\alpha _{t}$ is a function, $\alpha _{t}:\mathbb{Z}_{2}^{n}\times 
\mathbb{Z}_{2}^{n}\rightarrow \{-1,1\}\subset \mathbb{R}$, called the 
\textit{twist map}.

We remark that the elements of the group $\mathbb{Z}_{2}^{n}$ \ can be
considered as integers from $0$ to $2^{n}-1$ with multiplication "\textit{%
exclusive or}" for the binary representations. It is clear that this
operation is equivalent to the addition \ in $\mathbb{Z}_{2}^{n}$.\medskip

\textbf{Proposition 2.3.} \textit{Let} $\mathbb{N}=\{0,1,2,...\}$ \textit{be
the set of natural numbers and} $p,q,r,s,t\in \mathbb{N}$. \textit{We denote
with} $2p$ \textit{the twice of the number} $p,$ \textit{with} $\left(
p\right) ~$\textit{the binary representation of} $p$,\textit{\ with} $\left(
2p\right) \left( 2q\right) $ \textit{and} $pq$ \textit{the} "\textit{%
exclusive or"} \textit{for the binary representation of the numbers} $2p$ 
\textit{and} $2q$, \textit{respectively} $p$ \textit{and} $q$. \textit{%
Therefore, the following relations are true:}

i) $\left( 2p\right) \left( 2q\right) =2pq;$

ii) $\left( 2p\right) \left( 2q+1\right) =2pq+1;$

iii) $\left( 2p+1\right) \left( 2q\right) =2pq+1;$

iv) $\left( 2p+1\right) \left( 2q+1\right) =2pq;$

v) If $p,q\leq 2^{t}-1$, \textit{then} $p\left( 2^{t}+q\right) =2^{t}+pq;$

vi) If $p,r\leq 2^{t}-1$, \textit{then} $\left( 2^{t}+r\right) p=2^{t}+rp;$

vii) If $r,s\leq 2^{t}-1$, \textit{then} $pq=\left( 2^{t}+r\right) \left(
2^{t}+s\right) =rs.\medskip $

\textbf{Proof.} By straightforward calculation, since $\left( 2p\right) $ is
a moving of one position to the left of $\left( p\right) $, with $\left(
p\right) $ the binary representation of the number $p$. $_{{}}\medskip $

In the following, we consider the algebra $\mathcal{E}_{t}$ and arbitrary $%
\gamma _{1},...,\gamma _{t}\in K$, with $K$ a field of characteristic zero.
Let $B=\{f_{0},f_{1},...,f_{2^{t}-1}\}$ be the basis in the algebra $%
\mathcal{E}_{t}$. Using this basis, the basis in the algebra $\mathcal{E}%
_{t+1}$ can be written under the form:

$f_{0}=1=\left( 1,0\right) ,f_{1}=\left( f_{1},0\right)
,...,f_{2^{t}-1}=\left( f_{2^{t}-1},0\right) ,$

$f_{2^{t}}=\left( 0,1\right) ,f_{2^{t}+1}=\left( 0,f_{1}\right)
,f_{2^{t}+2}=\left( 0,f_{2}\right) ,...,$

$f_{2^{t}+i}=\left( 0,f_{i}\right) ,...,f_{2^{t+1}-1}=\left(
0,f_{2^{t}-1}\right) $.

We will prove that $\mathcal{E}_{t}$ is a twisted group algebra for the
group $G=\mathbb{Z}_{2}^{n},n=2^{t}$, $t\in \mathbb{N}$, over a field $K$,
with \textit{char}$K=0$.\medskip

\textbf{Example 2.4}. For $t=1$, we have the multiplication table

\begin{tabular}{c||c|c|}
$\cdot $ & $1$ & $\,\,\,f_{1}$ \\ \hline\hline
$\,1$ & $1$ & $\,\,\,f_{1}$ \\ \hline
$\,f_{1}$ & $\,\,f_{1}$ & $\gamma _{1}$ \\ \hline
\end{tabular}%
.

We remark that $\alpha _{1}\left( 0,0\right) =1,\alpha _{1}\left( 0,1\right)
=\alpha _{1}\left( 1,0\right) ,\alpha _{1}\left( 1,1\right) =\gamma _{1}$.
Moreover, $f_{p}f_{q}=\alpha _{1}\left( p,q\right) f_{p\otimes q}=f_{pq}$.
Indeed, since $1\otimes 0=0\otimes 1=1$ and $0\otimes 0=0$, we obtain that $%
\mathcal{E}_{t}$ is a twisted algebra.\medskip

\textbf{Example 2.5}. For $t=2$, we have the following multiplication table 
\begin{tabular}{c||c|c|c|c|}
$\cdot $ & $1$ & $\,\,\,f_{1}$ & $\,\,\,\,\,f_{2}$ & $\,\,\,\,f_{3}$ \\ 
\hline\hline
$\,1$ & $1$ & $\,\,\,f_{1}$ & $\,\,\,\,f_{2}$ & $\,\,\,\,f_{3}$ \\ \hline
$\,f_{1}$ & $\,\,f_{1}$ & $\gamma _{1}$ & $\,\,\,\,f_{3}$ & $\gamma
_{1}f_{2} $ \\ \hline
$\,f_{2}$ & $\,f_{2}$ & $-f_{3}$ & $\,\gamma _{2}$ & $\,\,-\gamma _{2}f_{1}$
\\ \hline
$f_{3}$ & $f_{3}$ & $-\gamma _{1}f_{2}$ & $\gamma _{2}f_{1}$ & $-\gamma
_{1}\gamma _{2}$ \\ \hline
\end{tabular}%
.

We remark that $\alpha _{2}\left( i,j\right) =-\alpha _{2}\left( j,i\right) $%
, for $\{1,2,3\}$ and $\alpha _{2}\left( 0,0\right) =\alpha _{2}\left(
0,i\right) =\alpha _{2}\left( i,0\right) =1=\alpha _{2}\left( 1,2\right) $,
for $i\in \{1,2,3\}$. Moreover, we have $\alpha _{2}\left( 1,1\right)
=\alpha _{2}\left( 1,3\right) =\gamma _{1},\alpha _{2}\left( 2,2\right)
=\alpha _{2}\left( 3,2\right) =\gamma _{2},\alpha _{2}\left( 3,3\right)
=-\gamma _{1}\gamma _{2}.\medskip $

\textbf{Proposition 2.6.} \textit{For all} $t\in \mathbb{N}$\textit{, we can
define a twist map} $\alpha _{t}:\mathbb{Z}_{2}^{n}\times \mathbb{Z}%
_{2}^{n}\rightarrow K,$ $n=2^{t},$ \textit{such that for} $f_{p},f_{q},$ 
\textit{two vectors of the basis} $B,$ \textit{with} $pq$ \textit{%
representing the "exclusive or" for the indices }$p$ \textit{and} $q$\textit{%
, we have that} $f_{p}f_{q}=\alpha _{t}\left( p,q\right) f_{pq}$. \textit{%
For the twist map} $\alpha _{t}$\textit{, the following relations are true:}

1) \textit{If} $p,q\in \{0,1,...,2^{t}-1\}$, \textit{then} $\alpha
_{t+1}\left( p,q\right) =\alpha _{t}\left( p,q\right) ;$

2) \textit{If} $p\in \{0,1,...,2^{t}-1\},q\in
\{2^{t},2^{t}+1,...,2^{t+1}-1\},q=2^{t}+r$, \textit{then} $\alpha
_{t+1}\left( p,q\right) =-\alpha _{t}\left( p,r\right) $, \textit{for} $%
r\neq p$, $p\neq 0,r\neq 0$.$~$\textit{If} $r\neq p$ \textit{and} $r=0$, 
\textit{then} $\alpha _{t+1}\left( p,q\right) =\alpha _{t}\left( p,0\right)
=1$. \textit{If} $r=p\neq 0$, \textit{then} $\alpha _{t+1}\left( p,q\right)
=\alpha _{t}\left( p,p\right) $.

3) \textit{If} $p\in \{2^{t},2^{t}+1,...,2^{t+1}-1\},q\in
\{0,1,...,2^{t}-1\},p=2^{t}+r$ \textit{and} $r\neq q$, \textit{then} $\alpha
_{t+1}\left( p,q\right) =-\alpha _{t}\left( r,q\right) $. \textit{If} $r=q$, 
\textit{then} $\alpha _{t+1}\left( p,q\right) =-\alpha _{t}\left( r,r\right) 
$.

4) \textit{If} $p\in \{2^{t},2^{t}+1,...,2^{t+1}-1\},q\in
\{2^{t},2^{t}+1,...,2^{t}-1\},$ \textit{with} $p=2^{t}+r$, $q=2^{t}+s$ 
\textit{and} $r\neq s,r\neq 0,s\neq 0$, \textit{then} $\alpha _{t+1}\left(
p,q\right) =\gamma _{t+1}\alpha _{t}\left( r,s\right) $. \textit{If} $%
r=s\neq 0$, \textit{then} $\alpha _{t+1}\left( p,p\right) =-\gamma
_{t+1}\alpha _{t}\left( r,r\right) $. If $r=0$, $s\not=0,$ \textit{then} $%
\alpha _{t+1}\left( p,q\right) =$ $-\gamma _{t+1}\alpha _{t}\left(
0,s\right) =-\gamma _{t+1}$. \textit{If} $s=0$, \textit{then} $\alpha
_{t+1}\left( p,q\right) =\gamma _{t+1}\alpha _{t}(r,0)=\gamma _{t+1}$%
.\medskip

\textbf{Proof}. By using induction over $t$. We assuming that the sentence
is true for $t$ and we will prove it for $t+1$.

\textbf{Case 1.} $p,q\in \{0,1,...,2^{t}-1\}$. In this situation, we have $%
f_{p}f_{q}=\left( f_{p},0\right) \left( f_{q},0\right) =\left(
f_{p}f_{q},0\right) =\left( \alpha _{t}\left( p,q\right) f_{pq},0\right)
=\alpha _{t}\left( p,q\right) \left( f_{pq},0\right) =\alpha _{t}\left(
p,q\right) f_{pq}$. Therefore $\alpha _{t+1}\left( p,q\right) =\alpha
_{t}\left( p,q\right) $.

\textbf{Case 2.} $p\in \{0,1,...,2^{t}-1\},q\in
\{2^{t},2^{t}+1,...,2^{t+1}-1\}$.

It results that $q=2^{t}+r$. Supposing that $r\neq p$, $p\neq 0,r\neq 0$, we
have $f_{p}f_{q}=f_{p}f_{2^{t}+r}=\left( f_{p},0\right) \left(
0,f_{r}\right) =$\newline
$=\left( 0,f_{r}f_{p}\right) =-\left( 0,f_{p}f_{r}\right) =-\left( 0,\alpha
_{t}\left( p,r\right) f_{pr}\right) =-\alpha _{t}\left( p,r\right) \left(
0,f_{pr}\right) =$\newline
$=-\alpha _{t}\left( p,r\right) f_{pq}$. We obtain $\alpha _{t+1}\left(
p,q\right) =-\alpha _{t}\left( p,r\right) $.

If $r\neq p$ and $r=0$, it results that $f_{p}f_{q}=f_{p}f_{2^{t}}=\left(
f_{p},0\right) \left( 0,f_{0}\right) =$\newline
$=\left( 0,f_{0}f_{p}\right) =\alpha _{t}\left( p,0\right) \left(
0,f_{p}\right) =\alpha _{t}\left( p,0\right) f_{2^{t}+p}$, therefore $\alpha
_{t+1}\left( p,q\right) =\alpha _{t}\left( p,0\right) $.

If $r=p\neq 0$, we have $f_{p}f_{q}=f_{p}f_{2^{t}+p}=\left( f_{p},0\right)
\left( 0,f_{p}\right) =\left( 0,f_{p}f_{p}\right) =\alpha _{t}\left(
p,p\right) \left( 0,1\right) =$\newline
$=\alpha _{t}\left( p,p\right) f_{2^{t}}$. Therefore, $\alpha _{t+1}\left(
p,q\right) =\alpha _{t}\left( p,p\right) $.

If $r\neq p$ and $r=0$, it results that $f_{p}f_{q}=f_{p}f_{2^{t}}=\left(
f_{p},0\right) \left( 0,f_{0}\right) =$\newline
$=\left( 0,f_{0}f_{p}\right) =\alpha _{t}\left( p,0\right) \left(
0,f_{p}\right) =\alpha _{t}\left( p,0\right) f_{2^{t}+p}$, therefore $\alpha
_{t+1}\left( p,q\right) =\alpha _{t}\left( p,0\right) $.

If $r=p\neq 0$, we have $f_{p}f_{q}=f_{p}f_{2^{t}+p}=\left( f_{p},0\right)
\left( 0,f_{p}\right) =\left( 0,f_{p}f_{p}\right) =\alpha _{t}\left(
p,p\right) \left( 0,1\right) =$\newline
$=\alpha _{t}\left( p,p\right) f_{2^{t}}$. Therefore, $\alpha _{t+1}\left(
p,q\right) =\alpha _{t}\left( p,p\right) $.

\textbf{Case 3.} $p\in \{2^{t},2^{t}+1,...,2^{t+1}-1\},q\in
\{0,1,...,2^{t}-1\}$. We have $p=2^{t}+r$. If $r\neq q$, therefore $%
f_{p}f_{q}=f_{2^{t}+r}f_{q}=\left( 0,f_{r}\right) \left( f_{q},0\right)
=\left( 0,f_{r}\overline{f_{q}}\right) =-\left( 0,f_{r}f_{q}\right) =-\left(
0,\alpha _{t}\left( r,q\right) f_{rq}\right) =$\newline
$=-\alpha _{t}\left( r,q\right) f_{pq}$. Therefore, $\alpha _{t+1}\left(
p,q\right) =-\alpha _{t}\left( r,q\right) $. If $r=q$, then $%
f_{p}f_{q}=f_{2^{t}+r}f_{r}=\left( 0,f_{r}\right) \left( f_{r},0\right)
=\left( 0,f_{r}\overline{f_{r}}\right) =-\alpha _{t}\left( r,r\right) \left(
0,1\right) =-\gamma _{r}f_{2^{t}}$. We obtain $\alpha _{t+1}\left(
p,q\right) =-\alpha _{t}\left( r,r\right) $.

\textbf{Case 4.} $p\in \{2^{t},2^{t}+1,...,2^{t+1}-1\},q\in
\{2^{t},2^{t}+1,...,2^{t+1}-1\},$ with $p=2^{t}+r$ and $q=2^{t}+s$. If $%
r\neq s$, $r\not=0,s\not=0,$we obtain $f_{p}f_{q}=f_{2^{t}+r}f_{2^{t}+s}=%
\left( 0,f_{r}\right) \left( 0,f_{s}\right) =\left( \gamma _{t+1}\overline{%
f_{s}}f_{r},0\right) =\gamma _{t+1}\left( f_{r}f_{s},0\right) =$\newline
$=\gamma _{t+1}\alpha _{t}\left( r,s\right) f_{rs}$. Therefore, $\alpha
_{t+1}\left( p,q\right) =\gamma _{t+1}\alpha _{t}(r,s)$.

If $r=s\neq 0$, then $f_{p}f_{q}=f_{2^{t}+r}f_{2^{t}+r}=\left(
0,f_{r}\right) \left( 0,f_{r}\right) =\left( \gamma _{t+1}\overline{f_{r}}%
f_{r},0\right) =-\gamma _{t+1}\left( f_{r}f_{r},0\right) =$\newline
$=-\gamma _{t+1}\alpha _{t}\left( r,r\right) f_{0}$. We get that $\alpha
_{t+1}\left( p,p\right) =-\gamma _{t+1}\alpha _{t}\left( r,r\right) $.

If $r=0$, it results $f_{p}f_{q}=f_{2^{t}}f_{2^{t}+s}=\left( 0,f_{0}\right)
\left( 0,f_{s}\right) =\left( \gamma _{t+1}\overline{f_{s}}f_{0},0\right) =$%
\newline
$=-\gamma _{t+1}\alpha _{t}\left( 0,s\right) \left( f_{s},0\right) =-\gamma
_{t+1}\alpha _{t}\left( 0,s\right) f_{s}$. Therefore, $\alpha _{t+1}\left(
p,q\right) =$ $-\gamma _{t+1}\alpha _{t}\left( 0,s\right) $.

If $s=0$, then $f_{p}f_{q}=f_{2^{t}+r}f_{2^{t}}=\left( 0,f_{r}\right) \left(
0,f_{0}\right) =\left( \gamma _{t+1}f_{0}f_{r},0\right) =$\newline
$=\gamma _{t+1}\alpha _{t}\left( r,0\right) f_{rs}$. Therefore, $\alpha
_{t+1}\left( p,q\right) =\gamma _{t+1}\alpha _{t}(r,0)$. If $r=s=0,\alpha
_{t+1}\left( p,p\right) =\alpha _{t+1}\left( 2^{t},2^{t}\right) $, then $%
f_{p}f_{q}=f_{2^{t}}f_{2^{t}}=\left( 0,f_{0}\right) \left( 0,f_{0}\right)
=\left( \gamma _{t+1}f_{0}f_{0},0\right) =$\newline
$=\gamma _{t+1}\alpha _{t}\left( 0,0\right) =\gamma _{t+1}._{{}}\medskip $

\textbf{Proposition 2.7.} \textit{For all} $t\in \mathbb{N}$\textit{, we can
define a sign map} $\theta _{t}:\mathbb{Z}_{2}^{n}\times \mathbb{Z}%
_{2}^{n}\rightarrow \{-1,1\},$ $n=2^{t},$ \textit{such that for} $%
f_{p},f_{q},$ \textit{two vectors of the basis} $B$\textit{, we have that} $%
f_{p}f_{q}$ \textit{has associated a sign,} $\theta _{t}\left( p,q\right) $. 
\textit{For the sign map} $\theta _{t}$, \textit{the following relations are
true:}

1) \textit{If} $p,q\in \{0,1,...,2^{t}-1\}$, \textit{then} $\theta
_{t+1}\left( p,q\right) =\theta _{t}\left( p,q\right) ;$

2) \textit{If} $p\in \{0,1,...,2^{t}-1\},q\in
\{2^{t},2^{t}+1,...,2^{t+1}-1\},q=2^{t}+r$, \textit{then} $\theta
_{t+1}\left( p,q\right) =-\theta _{t}\left( p,r\right) $, \textit{for} $%
r\neq p$, $p\neq 0,r\neq 0$.$~$\textit{If} $r\neq p$ \textit{and} $r=0$, 
\textit{then} $\theta _{t+1}\left( p,q\right) =\theta _{t}\left( p,0\right)
=+1$. \textit{If} $r=p\neq 0$, \textit{then} $\theta _{t+1}\left( p,q\right)
=\theta _{t}\left( p,p\right) $.

3) \textit{If} $p\in \{2^{t},2^{t}+1,...,2^{t+1}-1\},q\in
\{0,1,...,2^{t+1}-1\},p=2^{t}+r$ \textit{and} $r\neq q$, \textit{then} $%
\theta _{t+1}\left( p,q\right) =-\theta _{t}\left( r,q\right) $. \textit{If} 
$r=q$, \textit{then} $\theta _{t+1}\left( p,q\right) =-\theta _{t}\left(
r,r\right) $.

4) \textit{If} $p\in \{2^{t},2^{t}+1,...,2^{t+1}-1\},q\in
\{2^{t},2^{t}+1,...,2^{t+1}-1\}$, \textit{with} $p=2^{t}+r$, $q=2^{t}+s$ 
\textit{and} $r\neq s,r\neq 0,s\neq 0$, \textit{then} $\theta _{t+1}\left(
p,q\right) =\theta _{t}\left( r,s\right) $. \textit{If} $r=s\neq 0$, \textit{%
then} $\theta _{t+1}\left( p,p\right) =-\theta _{t}\left( r,r\right) $. 
\textit{If} $r=0$, $s\not=0,$ \textit{then} $\theta _{t+1}\left( p,q\right)
= $ $-\theta _{t}\left( 0,s\right) =-1$. \textit{If} $s=0$, \textit{then} $%
\theta _{t+1}\left( p,q\right) =\theta _{t}(r,0)=+1$.\medskip

\textbf{Proof.} It results from the above Proposition, by taking $\gamma
_{1}=\gamma _{2}=...=\gamma _{t+1}=1$.$_{{}}\medskip $

\textbf{Theorem 2.8.} \textit{The algebra} $\mathcal{E}_{t}$ \textit{is a
twisted group algebra for the group} $G=\mathbb{Z}_{2}^{n},n=2^{t}$\textit{,}
$t\in \mathbb{N}$, \textit{over a field} $K$\textit{,} \textit{with} \textit{%
char}$K=0$.\medskip \medskip

\textbf{Proof.} It results from Proposition 2.6 and Proposition 2.7.$%
_{{}}\medskip \medskip $

\textbf{Remark 2.9.} i) With the notations from the above propositions, let $%
\left( p\right) =i_{t}i_{t-1}...i_{1}$ and$~\left( q\right)
=j_{t}j_{t-1}...j_{1}$ be the binary representations for the indices of the
basis $B$. The coefficients $\alpha _{t}\left( p,q\right) $ has a sign,
obtained by using Proposition 2.7, and a product of elements from the set $%
\{1,\gamma _{1},\gamma _{2},...,\gamma _{t}\}$. From multiplication formula $%
\left( 1.1\right) ,$\ in algebra $\mathcal{E}_{t},$the element $\gamma _{m},$
$m\in \{1,2,...,t\},$ appears at the step $m$ and appears as coefficient of
the product $f_{p}f_{q}=\alpha _{t}\left( p,q\right) f_{pq}$ if and only if $%
i_{m}=j_{m}=1$. Indeed, if $i_{m}=j_{m}=1$, then, since the multiplication
at the step $m$ of the Cayley-Dickson process is $\left( a_{1},a_{2}\right)
\left( b_{1},b_{2}\right) =\left( a_{1}b_{1}+\gamma _{m}\overline{b}%
_{2}a_{2},a_{2}\overline{b_{1}}+b_{2}a_{1}\right) $, we have $%
f_{p}f_{q}=f_{2^{m}+r}f_{2^{m}+s}=\left( 0,f_{r}\right) \left(
0,f_{s}\right) =\left( \gamma _{m}\overline{f_{s}}f_{r},0\right) =\left(
\gamma _{m}f_{r}f_{s},0\right) =\gamma _{m}\alpha _{m-1}\left( r,s\right)
f_{rs}$, with $r,s\leq 2^{t}-1,r,s\in \mathbb{N}$.

ii) We remark that $\theta _{1}\left( i,j\right) =+1$, for $i,j\in \{0,1\}$%
.\medskip 

\textbf{The Algorithm\medskip }

Let $B=\{f_{0},f_{1},...,f_{2^{t}-1}\}$ be the basis in the algebra $%
\mathcal{E}_{t}$. We want compute $f_{p}f_{q},p,q\in \{0,1,2,...,2^{t}-1\}$.

Let $\left( p\right) =i_{t}i_{t-1}...i_{1}$ and$~\left( q\right)
=j_{t}j_{t-1}...j_{1}$ be the binary representations for the indices $p$ and 
$q$. Let $k_{1},k_{2},...,k_{z}$ be those indices such that $%
i_{k_{1}}=j_{k_{1}}=1,...,i_{k_{z}}=j_{k_{z}}=1.$

Therefore, $f_{p}f_{q}=\theta _{t}(p,q)\gamma _{k_{1}}\gamma
_{k_{2}}...\gamma _{k_{z}}f_{pq}$, with the sign $\theta _{t}(p,q)$ given by
the Proposition 2.6.\medskip\ 

\textbf{Example 2.10.} 

i) We consider the generalized octonion algebra.

1) We compute $f_{3}f_{5}$. The binary representations for $3$ and $5$ are $%
\left( 3\right) =011,\left( 5\right) =101$. Therefore we have $%
i_{1}=j_{1}=1\rightarrow \gamma _{1}$ and $\left( 3\right) \left( 5\right)
=110\rightarrow 6$. For the sign $\theta _{3}(3,5)$ we have, $\theta
_{3}(3,5)=-\theta _{2}(3,1)=+\theta _{1}(1,1)=+1$. It results that $%
f_{3}f_{5}=\gamma _{1}f_{6}$.

2) We compute $f_{6}f_{7}$. The binary representations for $6$ and $7$ are $%
\left( 6\right) =110,\left( 7\right) =111$. Therefore we have $%
i_{2}=j_{2}=1\rightarrow \gamma _{2}$, $i_{3}=j_{3}=1\rightarrow \gamma _{2}$
and $\left( 6\right) \left( 7\right) =001\rightarrow 1$. For the sign $%
\theta _{3}(6,7)$ we have, $\theta _{3}(6,7)=\theta _{2}(2,3)=-\theta
_{1}(0,1)=-1$. It results that $f_{6}f_{7}=-\gamma _{2}\gamma _{3}f_{1}$.

3) We compute $f_{6}f_{2}$. The binary representations for $6$ and $2$ are $%
\left( 6\right) =110,\left( 2\right) =010$. Therefore we have $%
i_{2}=j_{2}=1\rightarrow \gamma _{2}$ and $\left( 6\right) \left( 2\right)
=100\rightarrow 4$. For the sign $\theta _{3}(6,2)$ we have, $\theta
_{3}(6,2)=-\theta _{2}(2,2)=-\theta _{1}(1,1)=-1$. It results that $%
f_{6}f_{2}=-\gamma _{2}f_{4}$.

ii) We consider the generalized sedenion algebra.

1) We compute $f_{4}f_{12}$. The binary representations for $4$ and $12$ are$%
~\left( 4\right) =0100,\left( 12\right) =1100$. Therefore, we have $%
i_{3}=j_{3}=1\rightarrow \gamma _{3}$ and $\left( 4\right) \left( 7\right)
=1000\rightarrow 8$. For the sign, we get $\theta _{4}\left( 4,12\right)
=\theta _{3}\left( 4,4\right) =\theta _{2}\left( 0,0\right) =+1$. It results
that~ $f_{4}f_{12}=\gamma _{3}f_{8}$.

2) We compute $f_{10}f_{3}$. The binary decomposition for $10$ and $3$ are $%
\left( 10\right) =1010,\left( 3\right) =0011$. Therefore, we have $%
i_{2}=j_{2}=1\rightarrow \gamma _{2}$ and $\left( 10\right) \left( 3\right)
=1001\rightarrow 9$. For the sign, we get $\theta _{4}\left( 10,3\right)
=-\theta _{3}\left( 2,3\right) =-\theta _{2}\left( 2,3\right) =$\newline
$=+\theta _{1}\left( 0,1\right) =+1$. It results that $f_{10}f_{3}=\gamma
_{2}f_{9}$.

3) We compute $f_{9}f_{14}$. The binary decomposition for $9$ and $14$ are $%
\left( 9\right) =1001,\left( 14\right) =1110$. Therefore, we have $%
i_{4}=j_{4}=1\rightarrow \gamma _{4}$ and $\left( 9\right) \left( 14\right)
=0111\rightarrow 7$. For the sign, we get $\theta _{4}\left( 9,14\right)
=\theta _{3}\left( 1,6\right) =-\theta _{2}\left( 1,2\right) =-\theta
_{1}\left( 1,1\right) =-1$. It results that $f_{9}f_{14}=-\gamma _{4}f_{7}$.

The reader can consult Table 1, Table 2 and Table 3.

\begin{equation*}
\end{equation*}

\bigskip \qquad 
\begin{tabular}{c||c|c|c|c|}
$\cdot $ & $1$ & $\,\,\,f_{1}$ & $\,\,\,\,\,f_{2}$ & $\,\,\,\,f_{3}$ \\ 
\hline\hline
$\,1$ & $1$ & $\,\,\,f_{1}$ & $\,\,\,\,f_{2}$ & $\,\,\,\,f_{3}$ \\ \hline
$\,f_{1}$ & $\,\,f_{1}$ & $\gamma _{1}$ & $\,\,\,\,f_{3}$ & $\gamma _{1}f_{2}
$ \\ \hline
$\,f_{2}$ & $\,f_{2}$ & $-f_{3}$ & $\,\gamma _{2}$ & $\,\,-\gamma _{2}f_{1}$
\\ \hline
$f_{3}$ & $f_{3}$ & $-\gamma _{1}f_{2}$ & $\gamma _{2}f_{1}$ & $-\gamma
_{1}\gamma _{2}$ \\ \hline
\end{tabular}

\textbf{Table 1.} Multiplication table for the generalized quaternions

\bigskip 

~{\footnotesize $%
\begin{tabular}{c||c|c|c|c|c|c|c|c|}
$\cdot $ & $1$ & $\,\,\,f_{1}$ & $\,\,\,\,\,f_{2}$ & $\,\,\,\,f_{3}$ & $%
\,\,\,\,f_{4}$ & $\,\,\,\,\,\,f_{5}$ & $\,\,\,\,\,\,f_{6}$ & $%
\,\,\,\,\,\,\,f_{7}$ \\ \hline\hline
$\,1$ & $1$ & $\,\,\,f_{1}$ & $\,\,\,\,f_{2}$ & $\,\,\,\,f_{3}$ & $%
\,\,\,\,f_{4}$ & $\,\,\,\,\,\,f_{5}$ & $\,\,\,\,\,f_{6}$ & $%
\,\,\,\,\,\,\,f_{7}$ \\ \hline
$\,f_{1}$ & $\,\,f_{1}$ & $\gamma _{1}$ & $\,\,\,\,f_{3}$ & $\gamma _{1}f_{2}
$ & $\,\,\,\,f_{5}$ & $\gamma _{1}f_{4}$ & $-\,\,f_{7}$ & $\,\,-\gamma
_{1}f_{6}$ \\ \hline
$\,f_{2}$ & $\,f_{2}$ & $-f_{3}$ & $\gamma _{2}$ & $\,\,-\gamma _{2}f_{1}$ & 
$\,\,\,\,f_{6}$ & $\,\,\,\,\,f_{7}$ & $\gamma _{2}f_{4}$ & $\gamma _{2}f_{5}$
\\ \hline
$f_{3}$ & $f_{3}$ & -$\gamma _{1}f_{2}$ & $\gamma _{2}f_{1}$ & $-\gamma
_{1}\beta $ & $\,\,\,\,f_{7}$ & $\gamma _{1}f_{6}$ & $\,\,\,-\gamma _{2}f_{5}
$ & $-\gamma _{1}\gamma _{2}f_{4}$ \\ \hline
$f_{4}$ & $f_{4}$ & $-f_{5}$ & $-\,f_{6}$ & $-\,\,f_{7}$ & $\,\gamma _{3}$ & 
$\,\,-\,\gamma _{3}f_{1}$ & $\,\,-\gamma _{3}f_{2}$ & $\,\,\,\,-\,\gamma
_{3}f_{3}$ \\ \hline
$\,f_{5}$ & $\,f_{5}$ & -$\gamma _{1}f_{4}$ & $-\,f_{7}$ & -$\,\gamma
_{1}f_{6}$ & $\gamma _{3}f_{1}$ & -$\gamma _{1}\gamma _{3}$ & $\gamma
_{3}f_{3}$ & $\,\gamma _{1}\gamma _{3}f_{2}$ \\ \hline
$\,\,f_{6}$ & $\,\,f_{6}$ & $\,\,\,\,f_{7}$ & $\,\,-\gamma _{2}f_{4}$ & $%
\gamma _{2}f_{5}$ & $\gamma _{3}f_{2}$ & $\,\,\,-\gamma _{3}f_{3}$ & -$%
\gamma _{2}\gamma _{3}$ & $-\gamma _{2}\gamma _{3}f_{1}$ \\ \hline
$\,\,f_{7}$ & $\,\,f_{7}$ & $\gamma _{1}f_{6}$ & $\,-\gamma _{2}f_{5}$ & $%
\gamma _{1}\gamma _{2}f_{4}$ & $\gamma _{3}f_{3}$ & $-\gamma _{1}\gamma
_{3}f_{2}$ & $\gamma _{2}\gamma _{3}f_{1}$ & $\gamma _{1}\gamma _{2}\gamma
_{3}$ \\ \hline
\end{tabular}%
$}

\bigskip \textbf{Table 2.} Multiplication table for the generalized octonions

{\ \begin{figure}[htp]\center
	\epsfig{file=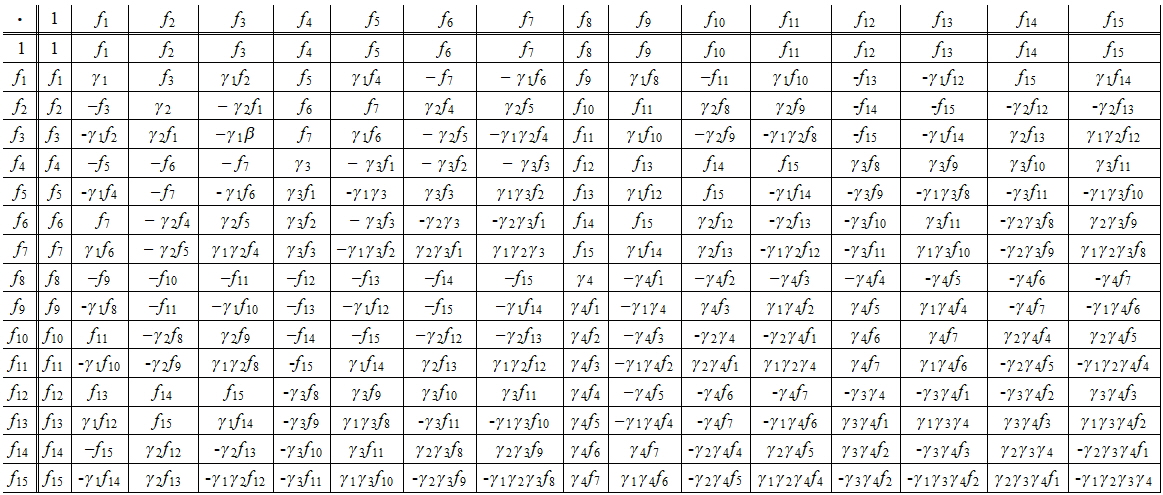,width=13cm}
\end{figure}
\newpage
\bigskip \textbf{Table 3.} Multiplication table for the generalized sedenions

\begin{equation*}
\end{equation*}

\textbf{3.} \textbf{Some remarks regarding nonassociative quaternion algebras%
}%
\begin{equation*}
\end{equation*}

In the following, we consider $K$ an arbitrary field with \textit{char}$%
K\not=2$. A nonassociative quaternion algebra is a $4$-dimensional $K$%
-algebra $A$ with identity, whose nucleus is equal to a separable quadratic
extension field $E$ of $K$. These algebras are division \ algebras, neither
quadratic, nor power-associative. Also, they are not third power-associative
(see [Pu; 14], [Wa; 87]).

Let $K\subset E$ \ be a separable quadratic field extension with $\sigma
:E\rightarrow E$, $\sigma \left( x\right) =\overline{x}$, an involution,
that means an automorphism which fixes $K$. We consider $\gamma \in E-K$. On
vector space \ $\mathbf{H}=E\oplus E$ we define the following multiplication 
\begin{equation}
\left( a_{1},a_{2}\right) \left( b_{1},b_{2}\right) =\left(
a_{1}b_{1}+\gamma (\overline{b}_{2}a_{2}),a_{2}\overline{b_{1}}%
+b_{2}a_{1}\right)   \tag{L}
\end{equation}%
Therefore $\mathbf{H}=E\oplus E$ becomes a nonassociative quaternion algebra
over $K$ with $\left( 1,0\right) $ as a unit element. (see [Pu; 14]) The
nucleus of this algebra is $E$.\medskip 

\textbf{Remark 3.1.} If we take $\alpha \in K$, let $E=K\left( \sqrt{\alpha }%
\right) $ and $\gamma \in E-K$. We have that $n\left( x\right) =x\sigma
\left( x\right) =x\overline{x},x\in E$. We remark that \thinspace $%
xx^{2}\neq x^{2}x$, for all $x\in \mathbf{H}$, then it is not (third)
power-associative, therefore is a nonassociative algebra and it is not a
flexible algebra.\medskip 

\textbf{Remark 3.2.} If, instead of $E$, we consider an arbitrary unitary
algebra $A$ over $K$ with the involution $\sigma :A\rightarrow A$ and $%
\gamma \in A$, an invertible element such that $\sigma \left( \gamma \right)
\not=\gamma $, then, on $\mathcal{A}=A\oplus A$, with a vector space
structure, we can get an algebra structure with the following multiplications

\begin{equation}
\left( a_{1},a_{2}\right) \left( b_{1},b_{2}\right) =\left(
a_{1}b_{1}+\gamma (\overline{b}_{2}a_{2}),a_{2}\overline{b_{1}}%
+b_{2}a_{1}\right)   \tag{L}
\end{equation}%
\begin{equation}
\left( a_{1},a_{2}\right) \left( b_{1},b_{2}\right) =\left( a_{1}b_{1}+%
\overline{b}_{2}(\gamma a_{2}),a_{2}\overline{b_{1}}+b_{2}a_{1}\right)  
\tag{M}
\end{equation}%
\begin{equation}
\left( a_{1},a_{2}\right) \left( b_{1},b_{2}\right) =\left( a_{1}b_{1}+(%
\overline{b}_{2}a_{2})\gamma ,a_{2}\overline{b_{1}}+b_{2}a_{1}\right)  
\tag{R}
\end{equation}%
depending on where is placed the element $\gamma $, used in the doubling
process. Therefore, we can continue this doubling process, obtaining
algebras with dimension double of dimension of $A$. These algebras are
denoted with $\mathcal{A}_{L}$, $\mathcal{A}_{M}$, $\mathcal{A}_{R}$,
depending on the chosen multiplication.These algebras are not isomorphic,
are nonassociative, are division and are not flexible algebras.(see [Pu;
14])\medskip\ 

\textbf{Proposition 3.3.} \textit{With the notations from Remark 3.1, let} 
\textbf{\ }$\mathbf{H}$ \textit{be a nonassociative quaternion algebra, with
basis} $\{1,f_{1},f_{2},f_{3}\}$\textit{, with involution} $\sigma $ \textit{%
such that} $f_{i}x=\sigma \left( x\right) f_{i},i\in \{1,2,3\}$, $\gamma \in
E-K$\textit{, with} $\sigma \left( \gamma \right) \not=\gamma ~$\textit{and
the multiplication given in the below table}

\begin{tabular}{c||c|c|c|c|}
$\cdot $ & $1$ & $\,\,\,f_{1}$ & $\,\,\,\,\,f_{2}$ & $\,\,\,\,f_{3}$ \\ 
\hline\hline
$\,1$ & $1$ & $\,\,\,f_{1}$ & $\,\,\,\,f_{2}$ & $\,\,\,\,f_{3}$ \\ \hline
$\,f_{1}$ & $\,\,f_{1}$ & $\alpha $ & $\,\,\,\,f_{3}$ & $\alpha f_{2}$ \\ 
\hline
$\,f_{2}$ & $\,f_{2}$ & $-f_{3}$ & $\,\gamma $ & $\,\,-\gamma f_{1}$ \\ 
\hline
$f_{3}$ & $f_{3}$ & $-\alpha f_{2}$ & $\gamma f_{1}$ & $-\alpha \gamma $ \\ 
\hline
\end{tabular}

\textit{We have that}

\textit{i)} $f_{i}^{2}f_{i}\not=f_{i}f_{i}^{2}$\textit{, for all} $i\in
\{2,3\}$;

\textit{ii) The elements of basis satisfy the following flexibility law}%
\begin{equation}
f_{i}\left( f_{k}f_{i}\right) =\left( f_{i}f_{k}\right) f_{i},\text{\textit{%
for all} }i,k\in \{1,2,3\},i\not=k.  \tag{F}
\end{equation}

\textbf{Proof.}

i)We have $f_{2}^{2}f_{2}=\gamma f_{2}$ and $f_{2}f_{2}^{2}=f_{2}\gamma
=\sigma \left( \gamma \right) f_{2}$. In the same way, we obtain that $%
f_{3}^{2}f_{3}=-\alpha \gamma f_{3}$ and $f_{3}f_{3}^{2}=-f_{3}\alpha \gamma
=-\alpha \sigma \left( \gamma \right) f_{3}$.

ii)Since $\gamma $ is in nuccleus of \textbf{\ }$\mathbf{H}$, we have:%
\newline
$f_{1}\left( f_{2}f_{1}\right) =f_{1}\left( -f_{3}\right) =-\alpha f_{2}$
and $\left( f_{1}f_{2}\right) f_{1}=f_{3}f_{1}=-\alpha f_{2}$;\newline
$f_{2}\left( f_{1}f_{2}\right) =f_{2}f_{3}=-\gamma f_{1}$ and $\left(
f_{2}f_{1}\right) f_{2}=-f_{3}f_{2}=-\gamma f_{1}$;\newline
$f_{1}\left( f_{3}f_{1}\right) =f_{1}\left( -\alpha f_{2}\right) =-\alpha
f_{3}$ and $(f_{1}f_{3})f_{1}=\alpha f_{2}f_{1}=-\alpha f_{3}\,$;\newline
$f_{3}\left( f_{1}f_{3}\right) =f_{3}\left( \alpha f_{2}\right) =\alpha
f_{3}f_{2}=\alpha \gamma f_{1}$ and $(f_{3}f_{1})f_{3}=-\alpha
f_{2}f_{3}=\alpha \gamma f_{1}$;\newline
$f_{2}\left( f_{3}f_{2}\right) =f_{2}\left( \gamma f_{1}\right) =\gamma
f_{2}f_{1}=-\gamma f_{1}f_{2}$ and $(f_{2}f_{3})f_{2}=\left( -\gamma
f_{1}\right) f_{2};$\newline
$f_{3}\left( f_{2}f_{3}\right) =-\gamma f_{3}f_{1}$ and $(f_{3}f_{2})f_{3}=%
\left( \gamma f_{1}\right) f_{3}=\gamma f_{1}f_{3}=-\gamma f_{3}f_{1}$.$%
_{{}}\medskip $

\textbf{Proposition 3.4.} \textit{Using notations from Remark 3.2, let} $A$ 
\textit{be an algebra of dimension} $n$ \textit{with basis} $%
\{1,f_{1},...,f_{n-1}\}$. \textit{If} $\gamma $ \textit{is in the nucleus of}
$A$ \textit{and the elements of basis satisfy the flexibility law,} $%
f_{i}\left( f_{k}f_{i}\right) =\left( f_{i}f_{k}\right) f_{i}$\textit{, for
all} $i,k\in \{1,2,...,n-1\},i\not=k$\textit{, therefore the elements of
basis of the algebra} $\mathcal{A}$ \textit{satisfy the same flexibility law 
}$\left( F\right) $\textit{.\medskip }

\textbf{Proof.} Since $\gamma $ is in the nucleus of $A$, therefore the
algebras $\mathcal{A}_{L},\mathcal{A}_{M},\mathcal{A}_{R}$ are the same. The
basis of the algebra $\mathcal{A}$ is $%
\{1,f_{1},...,f_{n-1},f_{n},...,f_{2n-1}\}$, where $f_{1}=\left(
f_{1},0\right) ,f_{2}=\left( f_{2},0\right) ,...,f_{n-1}=\left(
f_{n-1},0\right) ,f_{n}=\left( 0,1\right) ,f_{n+1}=\left( 0,f_{1}\right) ,$%
\newline
$f_{n+2}=\left( 0,f_{2}\right) ,...,f_{2n-1}=\left( 0,f_{n-1}\right) $.

\textbf{Case 1.} $i,k\in \{1,2,...,n-1\}$. We compute $f_{i}\left(
f_{k}f_{i}\right) $. We have $f_{i}\left( f_{k}f_{i}\right) =\left(
f_{i},0\right) (\left( f_{k},0\right) \left( f_{i},0\right) )=\left(
f_{i}\left( f_{k}f_{i}\right) ,0\right) =\left( (f_{i}f_{k})f_{i},0\right) =$%
\newline
$=(f_{i}f_{k})f_{i}$.

\textbf{Case 2.} $i\in \{1,2,...,n-1\}$ and $k\in
\{n,n+1,...,2n-1\},k=n+r,r>0$.

We obtain $f_{k}f_{i}=\left( 0,f_{r}\right) \left( f_{i},0\right) =\left(
0,f_{r}\overline{f_{i}}\right) =-\left( 0,f_{r}f_{i}\right) $ and $%
f_{i}\left( f_{k}f_{i}\right) =-\left( f_{i},0\right) \left(
0,f_{r}f_{i}\right) =$\newline
$=-\left( 0,\left( f_{r}f_{i}\right) f_{i}\right) =-\left( 0,-\left(
f_{r}f_{i}\right) f_{i}\right) =\left( 0,\left( f_{i}f_{r}\right)
f_{i}\right) $.

For the right term, we have $f_{i}f_{k}=\left( f_{i},0\right) \left(
0,f_{r}\right) =\left( 0,f_{r}f_{i}\right) $ and $\left( f_{i}f_{k}\right)
f_{i}=\left( 0,f_{r}f_{i}\right) \left( f_{i},0\right) =\left( 0,(f_{r}f_{i})%
\overline{f_{i}}\right) =-\left( 0,(f_{r}f_{i})f_{i}\right) =\left( 0,\left(
f_{i}f_{r}\right) f_{i}\right) $.

If $r=0$, we have $f_{i}\left( f_{k}f_{i}\right) =\left( f_{i},0\right)
\left( 0,f_{0}\overline{f_{i}}\right) =$\newline
$=\left( 0,\left( f_{0}\overline{f_{i}}\right) f_{i}\right) =\left(
0,-\left( f_{i}\right) f_{i}\right) =-\left( 0,f_{i}^{2}\right) $.

For the right term, we compute $f_{i}f_{k}=\left( f_{i},0\right) \left(
0,f_{r}\right) =\left( 0,f_{i}\right) $ and $\left( f_{i}f_{k}\right)
f_{i}=\left( 0,f_{i}\right) \left( f_{i},0\right) =\left( 0,f_{i}\overline{%
f_{i}}\right) =-\left( 0,f_{i}f_{i}\right) =-\left( 0,f_{i}^{2}\right) .$

\textbf{Case 3}. $k\in \{1,2,...,n-1\}$ and $i\in
\{n,n+1,...,2n-1\},i=n+s,s\geq 0$.

For $s\not=0$, we have $f_{k}f_{i}=\left( f_{k},0\right) \left(
0,f_{s}\right) =\left( 0,f_{s}f_{k}\right) $ and $f_{i}\left(
f_{k}f_{i}\right) =\left( 0,f_{s}\right) \left( 0,f_{s}f_{k}\right) =$%
\newline
$\left( \gamma \left( (\overline{f_{s}f\,_{k}})f_{s}\right) ,0\right)
=\left( \gamma \left( \left( f_{k}f_{s}\right) f_{s}\right) ,0\right)
=-\left( \gamma \left( \left( f_{s}f_{k}\right) f_{s}\right) ,0\right) $.

For the right term, we compute $f_{i}f_{k}=\left( 0,f_{s}\right) \left(
f_{k},0\right) _{r}=\left( 0,f_{s}\overline{f_{k}}\right) =-\left(
0,f_{s}f_{k}\right) $ and $\left( f_{i}f_{k}\right) f_{i}=-\left(
0,f_{s}f_{k}\right) \left( 0,f_{s}\right) =-\left( \gamma \left( \overline{%
f_{s}}\left( f_{s}f_{k}\right) \right) ,0\right) =-\left( \gamma \left(
f_{s}\left( f_{k}f_{s}\right) \right) ,0\right) =$\newline
$=-\left( \gamma \left( \left( f_{s}f_{k}\right) f_{s}\right) ,0\right) $,
because of flexibility of the elements $f_{s},f_{k},f_{s}$.

For $s=0$, we have $f_{k}f_{i}=\left( f_{k},0\right) \left( 0,1\right)
=\left( 0,f_{k}\right) $ and $f_{i}\left( f_{k}f_{i}\right) =\left(
0,1\right) \left( 0,f_{k}\right) =$\newline
$\left( \gamma \overline{f\,_{k}},0\right) =-\left( \gamma f_{k},0\right) $.

For the right term, we have $f_{i}f_{k}=\left( 0,1\right) \left(
f_{k},0\right) _{r}=\left( 0,\overline{f_{k}}\right) =-\left( 0,f_{k}\right) 
$ and $\left( f_{i}f_{k}\right) f_{i}=-\left( 0,f_{k}\right) \left(
0,1\right) =-\left( \gamma f_{k},0\right) .$

\textbf{Case 4.} $i\in \{n,n+1,...,2n-1\},k\in \{n,n+1,...,2n-1\}$. We have $%
k=n+r$ and $i=n+s,r,s>0$. It results $f_{k}f_{i}=\left( 0,f_{r}\right)
\left( 0,f_{s}\right) =\left( \gamma \left( \overline{f_{s}}f_{r}\right)
,0\right) ==-\left( \gamma \left( f_{s}f_{r}\right) ,0\right) $ and $%
f_{i}\left( f_{k}f_{i}\right) =-\left( 0,f_{s}\right) \left( \gamma \left(
f_{s}f_{r}\right) ,0\right) =\left( 0,f_{s}\left( \overline{\gamma \left(
f_{s}f_{r}\right) }\right) \right) =\left( 0,f_{s}((f_{r}f_{s})\overline{%
\gamma })\right) $. For the right term, we have $f_{i}f_{k}=\left(
0,f_{s}\right) \left( 0,f_{r}\right) =\left( \gamma \overline{f_{r}}%
f_{s},0\right) =-\left( \gamma f_{r}f_{s},0\right) $ and $\left(
f_{i}f_{k}\right) f_{i}=-\left( \gamma f_{r}f_{s},0\right) \left(
0,f_{s}\right) =-\left( 0,\gamma (f_{s}\left( f_{r}f_{s}\right) )\right)
=\left( 0,\gamma (\overline{f_{s}}\left( \overline{f_{r}}\overline{f_{s}}%
\right) )\right) =$\newline
$=\left( 0,\gamma (\overline{f_{s}}\left( \overline{f_{s}f_{r}}\right)
)\right) =\left( 0,\gamma \overline{(f_{s}f_{r})f_{s}}\right)
=(0,((f_{s}f_{r})f_{s})\overline{\gamma })$ and we have equality, due to the
flexibility of $f_{s},f_{r}$ and since $\gamma $ is in the nucleus of $A$.

If $s=0$, we have $f_{k}f_{i}=\left( 0,f_{r}\right) \left( 0,1\right)
=\left( \gamma f_{r},0\right) $ and $f_{i}\left( f_{k}f_{i}\right) =-\left(
0,1\right) \left( \gamma f_{r},0\right) =\left( 0,\left( \overline{\gamma
f_{r}}\right) \right) =-\left( 0,f_{r}\overline{\gamma }\right) =-\left(
0,\gamma f_{r}\right) $. For the right term, we have $f_{i}f_{k}=\left(
0,1\right) \left( 0,f_{r}\right) =\left( \gamma \overline{f_{r}},0\right)
=-\left( \gamma f_{r},0\right) $ and $\left( f_{i}f_{k}\right) f_{i}=-\left(
\gamma f_{r},0\right) \left( 0,1\right) =-\left( 0,\gamma f_{r}\right)
._{{}}\medskip $

\textbf{Remark 3.5.} From Proposition 3.3 and Proposition 3.4, if we
consider an element $\delta \in E$ and we apply the Cayley-Dickson process
for the algebra $\mathbf{H}$, we obtain a division octonion algebra in which
the elements of the basis satisfy the above flexibility law $\left( F\right) 
$.\medskip

\textbf{Example 3.6.} For $\mathbb{Q}\subset \mathbb{Q}\left( \sqrt{2}%
\right) $, we consider on $A=\mathbb{Q}\left( \sqrt{2}\right) \oplus \mathbb{%
Q}\left( \sqrt{2}\right) $ the multiplication given in the Proposition 3.3.
We have $i^{2}=2,j^{2}=\sqrt{2},k^{2}=-2\sqrt{2}$ and we get that $A$ over $%
\mathbb{Q}$ is a nonassociative quaternion division algebra, with basis $%
\{1,i,j,k\}$. Therefore, we obtain the following multiplication table:%
\newline
\newline

\begin{tabular}{c||c|c|c|c|}
$\cdot $ & $1$ & $\,\,\,i$ & $\,\,\,\,\,j$ & $\,\,\,\,k$ \\ \hline\hline
$\,1$ & $1$ & $\,\,\,i$ & $\,\,\,\,\,j$ & $\,\,\,\,k$ \\ \hline
$\,i$ & $\,i$ & $2$ & $\,\,\,\,k$ & $2j$ \\ \hline
$\,j$ & $\,j$ & $-k$ & $\,\sqrt{2}$ & $\,\,-i\sqrt{2}$ \\ \hline
$k$ & $k$ & $-2j$ & $i\sqrt{2}$ & $-2\sqrt{2}$ \\ \hline
\end{tabular}

We have:

$k=ij=\left( i,0\right) \left( 0,1\right) =\left( 0,i\right) $

$ik=\left( i,0\right) \left( 0,i\right) =\left( 0,i^{2}\right) =\left(
0,2\right) =2j$

$jk=\left( 0,1\right) \left( 0,i\right) =\left( -i\sqrt{2},0\right) =-i\sqrt{%
2}$

$k^{2}=\left( 0,i\right) \left( 0,i\right) =\left( \sqrt{2}i^{2},0\right) =-2%
\sqrt{2}$

We remark that \thinspace $jj^{2}\neq j^{2}j$. Indeed, since we have $%
jx=\sigma \left( x\right) j$, for all $x\in \mathbb{Q}\left( \sqrt{2}\right) 
$, it results \thinspace $jj^{2}=j\sqrt{2}=-\sqrt{2}j$ and $j^{2}j=\sqrt{2}j$%
, therefore are different. We compute $i\left( \sqrt{2}i\right) =-i\left( i%
\sqrt{2}\right) =i\left( jk\right) $ and $\left( i\sqrt{2}\right) i=-\left(
jk\right) i=i\left( jk\right) $. We denote this algebra with $A=\left( \frac{%
2,\sqrt{2}}{\mathbb{Q}}\right) $. Therefore, $A=\left( \frac{2,\sqrt{2}}{%
\mathbb{Q}}\right) $ is a nonassociative division quaternion algebra.

\begin{equation*}
\end{equation*}

\textbf{Conclusions. }In this paper we proved that an algebra obtained by
the Cayley-Dickson process is a twisted group algebra for the group $%
G=Z_{2}^{n},n=2^{t}$, $t\in N$, over a field $K$, with char$K=0$. Moreover,
we presented an algorithm which allows us to compose more easily two
elements from the basis, in this way the calculations become more easier in
higher dimension of the algebra $\mathcal{E}_{t}$. In the last section, we
give some properties and applications of the quaternion nonassociative
algebras. Since the quaternion nonassociatve algebras were not enough
exploited, we consider that the study of them can give us the chance to
obtain new and good results.

\begin{equation*}
\end{equation*}

\textbf{References}%
\begin{equation*}
\end{equation*}

[Ba; 09] Bales, J. W., \textit{A Tree for Computing the Cayley-Dickson Twist}%
, Missouri J. Math. Sci., \textbf{21(2)(}2009\textbf{)}, 83--93.

[Pu; 14] Pumpl\"{u}n S., \textit{How to obtain division algebras from a
generalized Cayley--Dickson doubling process}, Journal of Algebra,
402(2014), 406-434.

[Re; 71] Reynolds W. F., \textit{Twisted Group Algebras Over Arbitrary Fields%
}, Illinois Journal of Mathematics, 3(1971), 91-103.

[Sc; 66] Schafer R. D., \textit{An Introduction to Nonassociative Algebras,}
Academic Press, New-York, 1966.

[Sc; 54] Schafer R. D., \textit{On the algebras formed by the Cayley-Dickson
process,} Amer. J. Math., \textbf{76}(1954), 435-446.

[Wa; 87] Waterhouse W.C., \textit{Nonassociative quaternion algebras},
Algebras, Groups and Geometries, \textbf{4}(1987), 365 -- 378.

https://en.wikipedia.org/wiki/Change\_of\_rings%
\begin{equation*}
\end{equation*}

Cristina FLAUT

{\small Faculty of Mathematics and Computer Science, }\newline

{\small Ovidius University of Constan\c{t}a, Rom\^{a}nia,}

{\small Bd. Mamaia 124, 900527,}

{\small http://www.univ-ovidius.ro/math/}

{\small e-mail: cflaut@univ-ovidius.ro; cristina\_flaut@yahoo.com}

\begin{equation*}
\end{equation*}%
\qquad\ \ 

Remus BOBOESCU

{\small PhD student at Doctoral School of Mathematics,}

{\small Ovidius University of Constan\c{t}a, Rom\^{a}nia,}

{\small remus\_boboescu@yahoo.com}

\end{document}